\documentclass[11pt]{article}
\usepackage[centertags]{amsmath}
\usepackage{amsfonts}
\usepackage{amssymb}
\usepackage{amsthm}
\usepackage{amscd,amssymb,amsthm}
\newcommand{\be}{\begin{equation}}
\newcommand{\ee}{\end{equation}}
\newcommand{\bna}{\begin{eqnarray*}}
\newcommand{\ena}{\end{eqnarray*}}

\theoremstyle{definition}

\theoremstyle{remark}

\numberwithin{equation}{section}
\begin{document}
\title{ \bf   Radial distribution of Julia sets of derivatives of solutions of complex linear differential equations}
\author{\sc  Guowei Zhang$^1$, Jie Ding$^2$, Lianzhong Yang$^3$
 \\
1. School of Mathematics and Statistics, Anyang Normal University,
Anyang, Henan\\ 455000, P. R. China, E-mail: zhirobo@gmail.com\\
2. School of Mathematics, Taiyuan University of Technology, Taiyuan,
Shanxi \\ 030024, P. R. China, E-mail: jieding1984@gmail.com\\
3. School of Mathematics, Shandong University, Jinan, Shandong \\
250100, P. R. China, E-mail: lzyang@sdu.edu.cn
}
\date{}
\maketitle

\begin{center}{\bf Abstract}\end{center}

In this paper we mainly investigate the radial distribution of Julia
set of derivatives of entire solutions of some complex linear
differential equations. Under certain conditions, we find the lower
bound of it which improve some recent results.\\


\noindent{\bf 2000 MR subject classification}: 30D05; 37F10; 37F50

\baselineskip14pt \setcounter{section}{0}

\section{Introduction and main results}

In this paper, we assume the reader is familiar with standard
notations and basic results of Nevanlinna's value distribution
theory; see \cite{[Goldberg],Hayman,Laine,YL,YY}. Some basic
knowledge of complex dynamics of meromorphic functions is also
needed; see \cite{[Berg93],ZhengJH2}. Let $f$ be a meromorphic
function in the whole complex plane. We use $\sigma(f)$ and $\mu(f)$
to denote the order and lower order of $f$ respectively; see
\cite[p.10]{YY}
for the definitions.\\

We define $f^n, n\in \mathbb{N}$ denote the $n$th iterate of $f$.
The Fatou set $F(f)$ of  transcendental meromorphic function $f$ is
the subset of the plane $\mathbb{C}$ where the iterates $f^n$ of $f$
form a normal family. The complement of $F(f)$ in $\mathbb{C}$ is
called the Julia set $J(f)$ of $f$. It is well known that $F(f)$ is
open and completely invariant under $f$, $J(f)$ is closed and
non-empty.\\

We denote $\Omega(\alpha,\beta)=\{z\in\mathbb{C}|\arg
z\in(\alpha,\beta)\}$, where $0<\alpha<\beta<2\pi$. Given
$\theta\in[0,2\pi)$, if
$\Omega(\theta-\varepsilon,\theta+\varepsilon)\cap J(f)$ is
unbounded for any $\varepsilon>0$, then we call the ray $\arg
z=\theta$ the radial distribution of $J(f)$. Define
$$\Delta(f)=\{\theta\in[0,2\pi)|J(f) \rm {\ has\ the\ radial\ distribution\ with\ respect\ to} \arg z=\theta\}.$$

Obviously, $\Delta(f)$ is closed and so measurable. We use the
$meas\Delta(f)$ to denote the linear measure of $\Delta(f)$. Many
important results of radial distribution of transcendental
meromorphic functions have been obtained, for example
\cite{B3,Qiao1,Qiao2, QiuWu,WangS,ZhengJH3}.
 Qiao \cite{Qiao1}
proved that $meas\Delta(f)=2\pi$ if $\mu(f)<1/2$ and
$meas\Delta(f)\geq \pi/\mu(f)$ if $\mu(f)\geq1/2$, where $f(z)$ is a
transcendental entire function of finite lower order. Recently,
Huang et al \cite{Huang1, Huang2} considered radial distribution of
Julia set of entire solutions of linear complex
differential equations. Their results are stated as follows.\\

\noindent{\bf Theorem A} \cite{Huang1} {\it Let
$\{f_1,f_2,\ldots,f_n\}$ be a solution base of
\begin{eqnarray}\label{1.1} f^{(n)}+A(z)f=0,\end{eqnarray} where
$A(z)$ is a transcendental entire function with finite order, and
denote $E=f_1f_2\ldots f_n$.
Then $meas\Delta(E)\geq \min\{2\pi,\pi/\sigma(A)\}$}.\\

\noindent{\bf Theorem B} \cite{Huang2} {\it Let
$A_i(z)(i=0,1,\ldots,n-1)$ be entire functions of finite lower order
such that $A_0$ is transcendental and
$m(r,A_i)=o(m(r,A_0)),(i=1,2,\ldots,n-1)$ as $r\rightarrow\infty$.
Then every non-trivial solution $f$ of the equation
\begin{eqnarray}\label{1.2}
f^{(n)}+A_{n-1}f^{(n-1)}+\ldots+A_0f=0\end{eqnarray} satisfies
$meas\Delta(f)\geq \min\{2\pi,\pi/\mu(A_0)\}$.}
\\

 For entire functions and their derivatives, the difference between
 their local properties are astonishing, because a small disturbance
 of the parameter may cause a gigantic change of the dynamics for
 some given entire functions. So no one seems to believe
 that there are some neat relation between them in dynamical properties.
 However, Qiao \cite{Qiao3,Qiao2} proved that the Julia set of a transcendental entire function of
\emph{finite lower order} and its derivative have  a large amount of
common  radial distribution
 and their distribution densities influence each other.  A natural question is that what
 happens to the radial distribution of Julia set between entire function with  \emph{infinite lower order} and
 its derivative?
 \\

 It is easy to know that, by the logarithmic derivative lemma, the
non-trivial entire solutions of equations \eqref{1.1} and
\eqref{1.2} have infinite lower order, see details in \cite{Huang1}
and \cite{Huang2}. In the present paper, we study the radial
distribution of Julia set of the derivatives of entire solutions of
equations \eqref{1.1} and \eqref{1.2} and try to answer that above
question partially. Indeed,
we obtain the following results.\\

\noindent{\bf Theorem 1.1} {\it Let $A_i(z)(i=0,1,\ldots,n-1)$ be
entire functions of finite lower order such that $A_0$ is
transcendental and $m(r,A_i)=o(m(r,A_0)),(i=1,2,\ldots,n-1)$ as
$r\rightarrow\infty$. Then every non-trivial solution $f$ of the
equation \eqref{1.2} satisfies $meas(\Delta(f)\cap \Delta(f^{(k)}))
\geq \min\{2\pi,\pi/\mu(A_0)\}$,
where $k$ is a positive integer.}\\

\noindent{\bf Corollary 1.1} {\it Under the hypothesis of Theorem
1.1 we have $meas(\Delta(f^{(k)})) \geq \min\{2\pi,\pi/\mu(A_0)\}$,
where $k$ is a positive integer.}\\

Obviously, Theorem B is a corollary of Theorem 1.1. For entire solutions of equation \eqref{1.1}, we have\\

\noindent{\bf Corollary 1.2} {\it Assume that $f$ is any non-trivial
solution of equation \eqref{1.1}, we have $meas(\Delta(f^{(k)})) \geq
\min\{2\pi,\pi/\mu(A)\}$, where $k$ is a positive integer.}\\

Furthermore, we obtain the following.\\
\noindent{\bf Theorem 1.2} {\it Under the hypothesis of Theorem A,
we have $meas(\Delta(E^{(k)})) \geq \min\{2\pi,\pi/\sigma(A)\}$,
where $k$ is a positive integer.}\\

By Theorem 1.1, we have the next corollary even more.
\\

\noindent{\bf Corollary 1.3} {\it Suppose that
$A_i(z)(i=0,1,\ldots,n-1)$ be entire functions satisfying
$\sigma(A_j)<\mu(A_0)(j=1,2,\ldots,n-1)$ and $\mu(A_0)$ is finite.
Then every non-trivial solution $f$ of the equation \eqref{1.2}
satisfies $meas(\Delta(f)\cap \Delta(f^{(k)})) \geq
\min\{2\pi,\pi/\mu(A_0)\}$,
where $k$ is a positive integer.}\\

 \setcounter{section}{1}
\section{Preliminary lemmas}
At first, we recall the Nevanlinna characteristic in an angle; see
\cite{[Goldberg]}. We set
$$\Omega(\alpha,\beta,r)=\{z:z\in\Omega(\alpha,\beta),|z|<r\};$$
$$\Omega(r,\alpha,\beta)=\{z:z\in\Omega(\alpha,\beta),|z|\geq r\}$$
and denote by $\overline{\Omega}(\alpha,\beta)$ the closure of
$\Omega(\alpha,\beta)$. Let $g(z)$ be meromorphic on the angle
$\overline{\Omega}(\alpha,\beta)$, where $\beta-\alpha\in(0,2\pi]$.
Following \cite{[Goldberg]}, we define \begin{eqnarray}
A_{\alpha,\beta}(r,g)&=&\frac{w}{\pi}\int_1^r\left(\frac{1}{t^w}-\frac{t^w}{r^{2w}}\right)\{\log^+|g(te^{i\alpha})|+\log^+|g(te^{i\beta})|\}\frac{dt}{t};\nonumber\\
B_{\alpha,\beta}(r,g)&=&\frac{2w}{\pi
r^w}\int_{\alpha}^{\beta}\log^+|g(re^{i\theta})|\sin
w(\theta-\alpha)d\theta;\nonumber\\
C_{\alpha,\beta}(r,g)&=&2\sum_{1<|b_n|<r}\left(\frac{1}{|b_n|^w}-\frac{|b_n|^w}{r^{2w}}\right)\sin
w(\beta_n-\alpha),\nonumber \end{eqnarray} where
$w=\pi/(\beta-\alpha)$, and $b_n=|b_n|e^{i\beta_n}$ are poles of
$g(z)$ in $\overline{\Omega}(\alpha,\beta)$ appearing according to
their multiplicities. The Nevanlinna angular characteristic is
defined as
$$S_{\alpha,\beta}(r,g)=A_{\alpha,\beta}(r,g)+B_{\alpha,\beta}(r,g)+C_{\alpha,\beta}(r,g).$$
In particular, we denote the order of $S_{\alpha,\beta}(r,g)$ by
$$\sigma_{\alpha,\beta}(g)=\limsup_{r\rightarrow\infty}\frac{\log
S_{\alpha,\beta}(r,g)}{\log r}.$$\\

 We call $W$ is a hyperbolic domain if
$\overline{\mathbb{C}}\backslash W$ contains three points, where
$\overline{\mathbb{C}}$ is the extended complex plane. For an
$a\in\mathbb{C}\backslash W$, define
$$C_{W}(a)=\inf\{\lambda_W(z)|z-a|: \forall z\in W\},$$
where $\lambda_W(z)$ is the hyperbolic density on $W$. It's well
known that, if every component of $W$ is simply connected, then
$C_W(a)\geq 1/2$.\\

 \noindent{\bf Lemma 2.1.}\ (\cite[Lemma 2.2]{ZhengJH3}) {\it Let
$f(z)$ be an analytic in $\Omega(r_0,\theta_1,\theta_2)$, $U$ be a
hyperbolic domain, and $f:\Omega(r_0,\theta_1,\theta_2)\rightarrow
U$. If there exists a point $a\in \partial U\backslash \{\infty\}$
such that $C_U(a)>0$, then  there exists a constant $d>0$ such that,
for sufficiently small $\varepsilon>o$, we  have
$$|f(z)|=O(|z|^d),\ \ z\rightarrow\infty,\ z\in \Omega(r_0,\theta_1+\varepsilon,\theta_2-\varepsilon).
$$}

The next lemma shows some estimates for the logarithmic derivative
of functions being analytic in an angle. Before this, we recall the
definition of an R-set; for reference, see \cite{Laine}.  Set
$B(z_n,r_n)=\{z:|z-z_n|<r_n\}$. If $\sum_{n=1}^{\infty}r_n<\infty$
and $z_n\rightarrow\infty$, then $\cup_{n=1}^{\infty}B(z_n,r_n)$ is
called an R-set. Clearly, the set
$\{|z|:z\in\cup_{n=1}^{\infty}B(z_n,r_n)\}$ is of finite linear
measure.\\

\noindent{\bf Lemma 2.2.}\ (\cite[Lemma 2.2]{Huang2}) {\it Let
$z=re^{i\psi}, r_0+1<r$ and $\alpha\leq \psi\leq\beta$, where
$0<\beta-\alpha\leq2\pi$. Suppose that $n(\geq2)$ is an integer, and
that $g(z)$ is analytic in $\Omega(r_0,\alpha,\beta)$ with
$\sigma_{\alpha,\beta}(g)<\infty$. Choose
$\alpha<\alpha_1<\beta_1<\beta$. Then, for every
$\varepsilon_j\in(0,(\beta_j-\alpha_j)/2)(j=1,2,\ldots,n-1)$ outside
a set of linear measure zero with
$$\alpha_j=\alpha+\sum_{s=1}^{j-1}\varepsilon_s,\ \  \beta_j=\beta-\sum_{s=1}^{j-1}\varepsilon_s,\ \ j=2,3,\ldots,n-1.$$
there exists $K>0$ and $M>0$ only depending on $g$,
$\varepsilon_1,\ldots,\varepsilon_{n-1}$ and
$\Omega(\alpha_{n-1},\beta_{n-1})$, and not depending on $z$, such
that $$\left|\frac{g'(z)}{g(z)}\right|\leq Kr^M(\sin
k(\psi-\alpha))^{-2}$$ and
$$\left|\frac{g^{(n)}(z)}{g(z)}\right|\leq Kr^M\left(\sin
k(\psi-\alpha)\prod_{j=1}^{n-1}\sin
k_{\varepsilon_j}(\psi-\alpha_j)\right)^{-2}$$ for all $z\in
\Omega(\alpha_{n-1},\beta_{n-1})$ outside an R-set $D$, where
$k=\pi/(\beta-\alpha)$ and
$k_{\varepsilon_j}=\pi/(\beta_j-\alpha_j)(j=1,2,\ldots,n-1)$.}\\

\noindent{\bf Lemma 2.3.}  (\cite{YL2,ZhengJH2})\ {\it  Let $f(z)$
be a transcendental meromorphic function with lower order
$\mu(f)<\infty$ and order $0<\sigma(f)\leq\infty$. Then, for any
positive number $\lambda$ with $\mu(f)\leq \lambda\leq\sigma(f)$ and
any set $H$ of finite measure, there exists a sequence $\{r_n\}$
satisfies \\
(1). $r_n\not\in H, \lim_{n\rightarrow\infty}r_n/n=\infty$;\\
(2). $\liminf_{n\rightarrow\infty}\log T(r_n,f)/\log
r_n\geq\lambda$;\\
(3). $T(r,f)<(1+o(1))(2t/r_n)^{\lambda}T(r_n/2,f),
t\in[r_n/n,nr_n]$;\\
(4). $t^{-{\lambda-\varepsilon_n}}T(t,f)\leq
2^{\lambda+1}r_n^{-{\lambda-\varepsilon_n}}T(r_n,f), 1\leq t \leq
nr_n, \varepsilon_n=(\log n)^{-2}$.}\\

 Such $\{r_n\}$ is called a sequence of P\'{o}lya peaks of order
 $\lambda$ outside $H$. The following lemma, which related to P\'{o}lya
 peaks, is called the spread relation; see \cite{Baernstein}.\\

\noindent{\bf Lemma 2.4.}\ (\cite{Baernstein}){\it Let $f(z)$ be a
transcendental meromorphic function with positive order and finite
lower order, and has a deficient value $a\in \overline{\mathbb{C}}$.
Then, for any sequence of P\'{o}lya peaks $\{r_n\}$  of order
$\lambda>0,\ \mu(f)\leq \lambda\leq\sigma(f)$, and any positive
function$\Upsilon(r)\rightarrow0$ as $r_n\rightarrow\infty$, we have
$$\liminf_{r_n\rightarrow\infty}meas D_{\Upsilon}(r_n,a)\geq \min \left\{2\pi,\frac{4}{\lambda}\arcsin\sqrt{\frac{\delta(a,f)}{2}}\right\},$$
where
$$D_{\Upsilon}(r,a)=\left\{\theta\in[0,2\pi):\log^+\frac{1}{|f(re^{i\theta})-a|}>\Upsilon(r)T(r,f)\right\},\ a\in\mathbb{C}$$
 and
$$D_{\Upsilon}(r,\infty)=\left\{\theta\in[0,2\pi):\log^+|f(re^{i\theta})|>\Upsilon(r)T(r,f)\right\}.$$
 }

\setcounter{section}{2}
\section{Proof of Theorems}

\noindent{\bf Proof of theorem 1.1} We know that every non-trivial
solution $f$ of the equation is an entire function with infinite
lower order. We obtain the assertion by reduction to contradiction.
Assume that \begin{eqnarray}\label{3.1}
meas(\Delta(f)\cap\Delta(f^{(k)}))<\nu=\min\{2\pi,\pi/\mu(A_0)\}\end{eqnarray}
and so \begin{eqnarray}\label{3.2}
\xi:=\nu-meas(\Delta(f)\cap\Delta(f^{(k)}))>0.\end{eqnarray} Applying
Lemma 2.3 to $A_0$, we have a P\'{o}lya peak $\{r_j\}$ of order
$\mu(A_0)$ with all $r_j\not\in H$. Since $A_0$ is transcendental
entire function, it follows the Nevanlinna deficient
$\delta(\infty,A_0)=1$. By Lemma 2.4, for the P\'{o}lya peak
$\{r_j\}$, we have \begin{eqnarray}\label{3.3}
\liminf_{r_j\rightarrow\infty}meas(D_{\Upsilon}(r_j,\infty))\geq\pi/\mu(A_0),\end{eqnarray}
 where the function $\Upsilon(r)$ is defined by

 \begin{eqnarray}\label{3.4} \Upsilon(r)=\max \left\{\sqrt{\frac{\log r}{m(r,A_0)}},\sqrt{\frac{m(r,A_i)}{m(r,A_0)}},i=1,2,\ldots,n-1\right\}\end{eqnarray}
 and $m(r, A_j)$ is the proximation function of $A_j,
 j=0,1,\ldots,n-1$. Obviously, $\Upsilon(r)$ is positive and
 $\lim_{r\rightarrow\infty}\Upsilon(r)=0$. For the sake of
 simplicity, we denote $D_{\Upsilon}(r_j,\infty)$ by $D(r_j)$ in the
 following.
 We shall show that there must exist an open interval \begin{eqnarray}\label{3.5} I=(\alpha,\beta)\subset
\Delta(f^{(k)})^c,\ \ 0<\beta-\alpha<\nu\end{eqnarray} such that
\begin{eqnarray}\label{3.6} \lim_{j\rightarrow\infty} meas( \Delta(f)\cap
D(r_j)\cap I)>0, \end{eqnarray} where
$\Delta(f^{(k)})^c:=[0,2\pi)\backslash\Delta(f^{(k)})$. In order to
achieve this goal, we shall prove the following firstly.
\begin{eqnarray}\label{3.7}
\lim_{j\rightarrow\infty}meas(D(r_j)\backslash\Delta(f))=0.\end{eqnarray}
Otherwise, suppose that there is a subseries $\{r_{j_k}\}$ such that
\begin{eqnarray}\label{3.8} \lim_{k\rightarrow\infty}meas(D(r_{j_k})\backslash
\Delta(f))>0,\end{eqnarray} then there exists $\theta_0\in
\Delta(f)^c$ and $\eta>0$ satisfying \begin{eqnarray}\label{3.9}
\lim_{k\rightarrow\infty}meas((\theta_0-\eta,\theta_0+\eta)\cap
(D(r_{j_k})\backslash \Delta(f)))>0.\end{eqnarray} Since $\arg
z=\theta_0$ is not a radial distribution of $J(f)$, there exists
$r_0>0$ such that
\begin{eqnarray}\label{3.10} \Omega(r_0,\theta_0-\eta,\theta_0+\eta)\cap
J(f)=\emptyset.\end{eqnarray} This implies that there exists an
unbounded component $U$ of Fatou set $F(f)$, such that
$\Omega(r_0,\theta_0-\eta,\theta_0+\eta)\subset U$. Take a unbounded
and connected set $\Gamma\subset\partial U$, the mapping
$f:\Omega(r_0,\theta_0-\eta,\theta_0+\eta)\rightarrow
\mathbb{C}\backslash\Gamma$ is analytic. Since
$\mathbb{C}\backslash\Gamma$ is simply connected, then for any $a\in
\Gamma\backslash\{\infty\}$, we have
$C_{\mathbb{C}\backslash\Gamma}(a)\geq 1/2$. Now applying Lemma 2.1
to $f$ in $\Omega(r_0,\theta_0-\eta,\theta_0+\eta)$, for any
$\zeta>0,\zeta<\eta$, we have \begin{eqnarray}\label{3.11}
|f(z)|=O(|z|^{d_1}),\ \
z\in\Omega(r_0,\theta_0-\eta+\zeta,\theta_0+\eta-\zeta),\
|z|\rightarrow\infty,\end{eqnarray} where $d_1$ is a positive
constant. Recalling the definition of $S_{\alpha,\beta}(r,f)$, we
immediately get that \begin{eqnarray}\label{3.12}
S_{\theta_0-\eta+\zeta,\theta_0+\eta-\zeta}(r,f)=O(1).\end{eqnarray}
Therefore, by Lemma 2.2, there exists constants $M>0$ and $K>0$ such
that \begin{eqnarray}\label{3.13}
\left|\frac{f^{(s)}(z)}{f(z)}\right|\leq Kr^M, \ \
(s=1,2,\ldots,n-1),\end{eqnarray} for all
$z\in\Omega(r_0,\theta_0-\eta+\zeta,\theta_0+\eta-\zeta)$, outside a
R-set $H$.

Since $\zeta$ can be chosen sufficiently small, from \eqref{3.9} we
have \begin{eqnarray}\label{3.14}
\lim_{k\rightarrow\infty}meas((\theta_0-\eta+\zeta,\theta_0+\eta-\zeta)\cap
D(r_j))>0.\end{eqnarray} Thus, we can find an infinite series
$\{r_{j_k}e^{i\theta_{j_k}}\}$ such that for all sufficiently large
$k$, \begin{eqnarray}\label{3.15}
\log^+|A_0(r_{j_k}e^{i\theta_{j_k}})|>\Upsilon(r_{j_k})T(r_{j_k},A_0)=\Upsilon(r_{j_k})m(r_{j_k},A_0)
\end{eqnarray} where $\theta_{j_k}\in
F_{j_k}:=(\theta_0-\eta+\zeta,\theta_0+\eta-\zeta)\cap D(r_j)$.
Then, for sufficiently large $k$, we have
\begin{eqnarray}\label{3.16}
\int_{F_{j_k}}\log^+|A_0(r_{j_k}e^{i\theta_{j_k}})|d\theta\geq meas
(F_{j_k})\Upsilon(r_{j_k})m(r_{j_k},A_0). \end{eqnarray} On the
other hand, combining \eqref{1.2} and
 \eqref{3.9} leads to \begin{eqnarray}\label{3.17}
\int_{F_{j_k}}\log^+|A_0(r_{j_k}e^{i\theta_{j_k}})|d\theta&\leq&
\int_{F_{j_k}}\left(\sum_{s=1}^n
\log^+\left|\frac{f^{(s)}(r_{j_k}e^{i\theta_{j_k}})}{f(r_{j_k}e^{i\theta_{j_k}})}\right|+\sum_{i=1}^{n-1}\log^+|A_i(r_{j_k}e^{i\theta_{j_k}})|\right)d\theta+O(1)\nonumber\\
&=&
\int_{F_{j_k}}\left(\sum_{i=1}^{n-1}\log^+|A_i(r_{j_k}e^{i\theta_{j_k}})|\right)d\theta+O(\log
r_{j_k})\nonumber\\
&\leq&\sum_{i=1}^{n-1}m(r_{j_k},A_i)+O(\log r_{j_k})\nonumber\\
&\leq&c_0\Upsilon^2(r_{j_k})m(r_{j_k},A_0)\end{eqnarray} where $c_0$
is a positive constant. From \eqref{3.16} and \eqref{3.17}, we have
\begin{eqnarray}\label{3.18} 0<meas (F_{j_k})\leq c_0\Upsilon(r_{j_k})\end{eqnarray} which
contradicts to the fact $\Upsilon(r_{j_k})\rightarrow0$ as
$k\rightarrow\infty$. This contradiction implies \eqref{3.7} is
valid. By Theorem B, we know that \begin{eqnarray}\label{3.19} meas
\Delta(f)\geq \nu.\end{eqnarray}  From Lemma 2.4, we have, for all
sufficiently large $j$ and any positive $\varepsilon$,
\begin{eqnarray}\label{3.20} meas
D(r_j)>\nu-\varepsilon.\end{eqnarray} Combining \eqref{3.7},
\eqref{3.19} and \eqref{3.20} follows that, for all sufficiently
large $j$,
\begin{eqnarray}\label{3.21}  meas (\Delta(f)\cap D(r_j))\geq \nu-\xi/4,\end{eqnarray}
where $\xi$ is defined in \eqref{3.2}. Since $\Delta(f^{(k)})$ is
closed, clearly $\Delta(f^{(k)})^c$ is open, so it consists of at
most countably open intervals. We can choose finitely many open
intervals $I_j,(j=1,2,\ldots,m)$ satisfying
\begin{eqnarray}\label{3.22} I_j\subset \Delta(f^{(k)})^c,\ \
meas(\Delta(f^{(k)})^c\backslash\cup_{i=1}^m
I_i)<\xi/4.\end{eqnarray}
 Since, for sufficiently large
$j$,  \begin{eqnarray}\label{3.23} & &meas( \Delta(f)\cap
D(r_j)\cap(\cup_{i=1}^m
I_i))+meas(\Delta(f)\cap D(r_j)\cap \Delta(f^{(k)}))\nonumber\\
&=&meas(\Delta(f)\cap D(r_j)\cap(\Delta(f^{(k)})\cup(\cup_{i=1}^m
I_i)))\geq \nu-\xi/2, \end{eqnarray} we have
\begin{eqnarray}\label{3.24} meas( \Delta(f)\cap
D(r_j)\cap(\cup_{i=1}^m I_i))&\geq&
\nu-\xi/2-meas(\Delta(f)\cap D(r_j)\cap \Delta(f^{(k)}))\nonumber\\
&\geq& \nu-\xi/2-meas(\Delta(f)\cap \Delta(f^{(k)}))=\xi/2.
\end{eqnarray} Thus, there exists an open interval
$I_{i_0}=(\alpha,\beta)\subset \cup_{i=1}^m I_i\subset
\Delta(f^{(k)})^c$ such that, for infinitely many sufficiently large
$j$, \begin{eqnarray}\label{3.25} meas( \Delta(f)\cap D(r_j)\cap
I_{i_0})\geq \frac{\xi}{2m}>0.\end{eqnarray} Then, we prove
\eqref{3.6} holds.\\

From \eqref{3.6}, we know that there are $\widetilde{\theta_0}$ and
$\widetilde{\eta}>0$ such that
\begin{eqnarray}\label{3.26}(\widetilde{\theta_0}-\widetilde{\eta},\widetilde{\theta_0}+\widetilde{\eta})\subset
I\end{eqnarray} and \begin{eqnarray}\label{3.27}
\lim_{j\rightarrow\infty} meas( \Delta(f)\cap D(r_j)\cap
(\widetilde{\theta_0}-\widetilde{\eta},\widetilde{\theta_0}+\widetilde{\eta}))>0.\end{eqnarray}
Then, there exists $\widetilde{r_0}$ such that
$\Omega(\widetilde{r_0},\widetilde{\theta_0}-\widetilde{\eta},\widetilde{\theta_0}+\widetilde{\eta})\cap
J(f^{(k)}(z))=\emptyset$. By the similar argument between
\eqref{3.10} and \eqref{3.11}, for any
$\widetilde{\zeta}>0,\widetilde{\zeta}<\widetilde{\eta}$, we have
\begin{eqnarray}\label{3.28} |f^{(k)}(z)|=O(|z|^{d_2}),\ \
z\in\Omega(\widetilde{r_0},\widetilde{\theta_0}-\widetilde{\eta}+\widetilde{\zeta},\widetilde{\theta_0}+\widetilde{\eta}-\widetilde{\zeta}),\
|z|\rightarrow\infty,\end{eqnarray} where $d_2$ is a positive
constant. By \eqref{3.27} we can choose an unbounded series
$\{r_je^{i\theta_j}\}$, for all sufficiently large $j$ such that
\begin{eqnarray}\label{3.29}
\log^+|A_0(r_je^{i\theta_j})|>\Upsilon(r_j)m(r_j,A_0),\end{eqnarray}
where
$$\theta_j\in  \Delta(f)\cap
D(r_j)\cap
(\widetilde{\theta_0}-\widetilde{\eta},\widetilde{\theta_0}+\widetilde{\eta}).$$
Fixed $r_Je^{i\theta_J}$, and  take a
$r_je^{i\theta_j}\in\{r_je^{i\theta_j}\}$. Take a simple Jordan arc
$\gamma$ in
$\Omega(\widetilde{r_0},\widetilde{\theta_0}-\widetilde{\eta},\widetilde{\theta_0}+\widetilde{\eta})$
which connecting $r_Je^{i\theta_J}$ to $r_Je^{i\theta_j}$ along
$|z|=r_J$, and connecting $r_Je^{i\theta_j}$ to $r_je^{i\theta_j}$
along $\arg z=\theta_j$. For any $z\in\gamma$, $\gamma_z$ denotes a
part of $\gamma$, which connecting $r_Je^{i\theta_J}$ to $z$. Let
$L(\gamma)$ be the length of $\gamma$. Clearly,
$$L(\gamma)=O(r_j),\ \ j\rightarrow\infty.$$
By \eqref{3.28}, it follows \begin{eqnarray}|f^{(k-1)}(z)|&\leq&
\int_{\gamma_z}|f^{(k)}(z)||dz|+c_k\nonumber \\
&\leq&O(|z|^{d_2}L(\gamma))+c_k\nonumber\\&\leq&O(r_j^{d_2+1}), \ \
j\rightarrow\infty.\nonumber\end{eqnarray} Similarly, we have
\begin{eqnarray}\label{3.30}|f^{(k-2)}(z)|&\leq&
\int_{\gamma_z}|f^{(k-1)}(z)||dz|+c_{k-1}\nonumber \\
&\leq&O(r_j^{d_2+2}), \ \ j\rightarrow\infty.\nonumber\\ &\vdots&\nonumber\\
|f(z)|&\leq&
\int_{\gamma_z}|f'(z)||dz|+c_1\nonumber \\
&\leq&O(r_j^{d_2+k}), \ \ j\rightarrow\infty.\end{eqnarray} where
$c_1,c_2,\ldots,c_k$ are constants, which are independent of $j$.
Therefore,\begin{eqnarray}\label{3.31}
S_{\widetilde{\theta_0}-\widetilde{\eta}+\widetilde{\zeta},\widetilde{\theta_0}+\widetilde{\eta}-\widetilde{\zeta}}(r,f)=O(1).\end{eqnarray}
By Lemma 2.2, we know \eqref{3.13} also holds for  all
$z\in\Omega(\widetilde{r_0},\widetilde{\theta_0}-\widetilde{\eta}+\widetilde{\zeta},\widetilde{\theta_0}+\widetilde{\eta}-\widetilde{\zeta})$,
outside a R-set $H$. Combining \eqref{3.13} and \eqref{3.29}, and
applying the similar argument as \eqref{3.16} and \eqref{3.17}, we
can deduce a contradiction. Therefore, it follows
\begin{eqnarray}\label{3.33}
meas(\Delta(f)\cap\Delta(f^{(k)}))\geq\min\{2\pi,\pi/\mu(A)\}.\end{eqnarray}
The proof is complete.\\

\noindent{\bf Proof of theorem 1.2}\ \   The main idea of this proof
comes from that of the proof of Theorem 1.1 in \cite{Huang1}, but
need some changes. We assume that
$meas(\Delta(E^{(k)}))<\min\{2\pi,\pi/\sigma(A)\}$. By similar
argument in \cite{Huang1}, there exists an angular domain
$\Omega(\alpha,\beta)$ such that \begin{eqnarray}
\Omega(\alpha,\beta)\cap \Delta(E^{(k)})=\emptyset,\ \
\Omega(r_0,\alpha,\beta)\cap J(E^{(k)})=\emptyset \end{eqnarray} for
sufficiently large $r_0$. Then by the same method between
\eqref{3.10} and \eqref{3.11}, we have \begin{eqnarray}
|E^{(k)}(z)|=O(|z|^{d}),\ \ z\in\Omega(r_0,\alpha,\beta),\
|z|\rightarrow\infty,\end{eqnarray} where $d$ is a positive
constant. Take a simple Jordan arc $\gamma$, which connected points
$z_0$ and $z$, satisfying $\gamma\subset \Omega(r_0,\alpha,\beta)$.
Applying the method which is used in \eqref{3.30}, we obtain
\begin{eqnarray} |E(z)|=O(|z|^{d+k}),\ \
z\in\Omega(r_0,\alpha,\beta),\ |z|\rightarrow\infty.\end{eqnarray}
Therefore, Theorem 1.2 can be proved word by word following the
proof of Theorem 1.1 in \cite{Huang1}.


\begin{thebibliography}{16}

\bibitem{Baernstein} A. Baernstein, Proof of Edrei's spread conjecture,
Proc. Lond. Math. Soc. 26(1973) 418-434.
\bibitem{B3} I. N. Baker, Sets of non-normality in
iteration theory, J. Lond. Math. Soc. 40(1965) 499-502.
\bibitem{[Berg93]}{}
W. Bergweiler, Iteration of meromorphic functions, Bull. Amer. Math.
Soc. (N. S.) 29(1993) 151-188.
\bibitem{[Goldberg]}{}
A. A. Gol'dberg, I. V. Ostrovskii, Value distribution of meromorphic
function, in: AMS Translations of Mathematical Monographs series,
2008.




\bibitem{Hayman}W. Hayman, Meromorphic functions, Clarendon Press, Oxford 1964.
\bibitem{Huang1}Z. Huang, J. Wang, On the radial distribution of Julia sets of entire solutions of $f^{(n)}+A(z)f=0$, J. Math. Anal. Appl. 387 (2012) 1106-1113.
\bibitem{Huang2}Z. Huang, J. Wang, On limit directions of Julia sets of entire solutions of linear differential equations, J. Math. Anal. Appl. 409 (2014) 478-484.



\bibitem{Laine}I. Laine, Nevanlinna Theory and Complex Differential Equations, de Gruyter,
Berlin, 1993.
%




%
%
%

\bibitem{Qiao3}J. Qiao, Julia set of entire functions and their derivatives, Chin. Sci. Bull. 39(3)(1994) 186-188.
\bibitem{Qiao1}J. Qiao, Stable domains in the iteration of entire functions(in Chinese), Acta. Math. Sin. 37 (1994) 702-708.
\bibitem{Qiao2}J. Qiao, On limiting directions of Julia set, Ann. Acad. Sci. Fenn. Math. 26 (2001) 391-399.
\bibitem{QiuWu}L. Qiu, S. J. Wu, Radial distributions of Julia sets of meromorphic functions, J. Aust. Math. Soc. 81 (3)(2006) 363-368.
\bibitem{WangS}S. Wang, On radial distributions of Julia sets of meromorphic functions, Taiwanese J. Math. 11(5) (2007) 1301-1313.
\bibitem{YL}L. Yang, Value distribution theory, Springer-Verlag, Berlin, 1993.
\bibitem{YL2}L. Yang, Borel directions of meromorphic functions in an
angular domain, Sci. China Math. Ser.(1)(1979) 149-163.
\bibitem{YY} C. C. Yang, H. X. Yi, Uniqueness theory of meromorphic functions, Science Press/Kluwer Academic Publishers, Beijing/Dordrecht,
2003.
\bibitem{ZhengJH}J. H. Zheng, On multiply-connected Fatou components in iteration of meromorphic functions,
J. Math. Anal. Appl., 313 (2006) 24-37
\bibitem{ZhengJH2}
J. H. Zheng, Dynamics of meromorphic functions (in Chinese),
Tsinghua University Press, Beijing, 2006.
\bibitem{ZhengJH3}
J. H. Zheng, S. Wang, Z. Huang, Some properties of Fatou and Julia
sets of transcendental meromorphic functions, Bull. Aust. Math. Soc.
66 (2002) 1-8.
\end{thebibliography}
\end{document}